\title{The M\"{o}bius Function of the Suzuki Groups,\\
with Applications to Enumeration}

\author{Martin Downs and Gareth A. Jones\\
School of Mathematics\\
University of Southampton\\
Southampton SO17  1BJ, U.K.\\
{\tt G.A.Jones@maths.soton.ac.uk}
}

\documentclass[12pt]{article}
\usepackage{color}
\usepackage[latin1]{inputenc}
\usepackage{a4wide}
\usepackage{latexsym}
\usepackage{amsmath}
\usepackage{amssymb}
\usepackage{amsfonts}
\newtheorem{thm}{Theorem}[section]
\newtheorem{lemma}[thm]{Lemma}

\newtheorem{prop}[thm]{Proposition}
\date{}
\begin{document}

\maketitle

\begin{abstract}
We compute the M\"obius function for the subgroup lattice of the simple Suzuki group $Sz(q)$; this is used to enumerate normal subgroups of certain Hecke groups with quotients isomorphic to $Sz(q)$. 
\end{abstract}

\noindent{\bf MSC classification:} 20G41(primary); 
05A99, 
11A25, 
14H57, 
20B20, 
20H10 
(secondary).

\medskip

\noindent{\bf Key words:} Suzuki group, M\"obius function, subgroup lattice, Hecke group.

\section{Introduction}\label{Intro}

Hall's theory of M\"obius inversion in groups~\cite{Hal} allows one to enumerate various objects associated with a given finite group $G$. For example, using this method one can compute the number of normal subgroups $N$ of the free group $F_k$ of finite rank $k$ with $F_k/N\cong G$, or equivalently the number of orbits of ${\rm Aut}\,G$ on generating $k$-tuples for $G$. Indeed, both of these numbers are equal to
\begin{equation}\label{d_k(G)}
d_k(G)=\frac{1}{|{\rm Aut}\,G|}\sum_{H\le G}\mu_G(H)|H|^k
\end{equation}
where $\mu_G$ is the M\"obius function on the lattice of subgroups of $G$, defined recursively by
\begin{equation}\label{mu}
\sum_{K\ge H}\mu_G(K)=\delta_{H,G}.
\end{equation}
Here $\delta_{H,G}$ is the Kronecker delta function, equal to $1$ or $0$ as $H=G$ or $H<G$.) As examples of this, $d_2(G)$ gives
\begin{itemize}
\item the number of isomorphism classes of orientably regular hypermaps with automorphism group $G$;
\item the number of regular unbranched coverings of the sphere minus three points (or of the torus minus one point) with covering group $G$;
\item the number of regular dessins (in Grothendieck's terminology~\cite{Gro}) with automorphism group $G$.
\end{itemize}

Equation~(\ref{d_k(G)}) is a particular case of the general result that the number $n_{\Gamma}(G)$ of normal subgroups $N$ of a finitely generated group $\Gamma$ with quotient isomorphic to a given finite group $G$ is given by
\begin{equation}\label{NGammaG}
n_{\Gamma}(G)
=\frac{|{\rm Epi}(\Gamma,G)|}{|{\rm Aut}\,G|}=\frac{1}{|{\rm Aut}\,G|}\sum_{H\le G}\mu_G(H)|{\rm Hom}(\Gamma,H)|,
\end{equation}
which arises from applying M\"obius inversion to the obvious equation
\begin{equation}
|{\rm Hom}(\Gamma,G)|=\sum_{H\le G}|{\rm Epi}(\Gamma,H)|.
\end{equation}
A similar principle applies to the enumeration of torsion-free normal subgroups of $\Gamma$ with quotient $G$: one simply counts the smooth homomorphisms and epimorphisms $\Gamma\to H$, those preserving the orders of torsion elements.

Implementing equation~(\ref{NGammaG}), as for some other similar results, for a specific group $G$ depends on knowing the value of $\mu_G(H)$ for each subgroup $H\le G$. Finding these values can be a laborious task, but once it has been done, equation~(\ref{NGammaG}) can be applied to $G$ in many different contexts, depending on the choice of $\Gamma$ (see~\cite{DJ, Jon95}, for example). The function $\mu_G$ has been determined for certain groups $G$, including the simple groups $L_2(p)=PSL_2(p)$ for primes $p\ge 5$, by Hall in~\cite{Hal}. Subsequently the first author extended this in~\cite{DowPhD} to the groups $L_2(q)$ and $PGL_2(q)$ for all prime powers $q$: see~\cite{DowJLMS} for full details for $L_2(2^e)$ and a statement of results for $L_2(q)$ with $q$ odd, and~\cite{DJ} for some applications.

The objective of this paper, which is based on an earlier preprint~\cite{Dow93}  by the first author, is to calculate the M\"obius function $\mu_G$ for the family of simple groups $G=Sz(q)={}^2B_2(q)$, where $q=2^e$ for some odd $e>1$. These groups were discovered in 1960 by Suzuki~\cite{Suz60, Suz}. We describe $G$ and its subgroups $H$ in Section~2, and give the values of $\mu_G(H)$ in Table~\ref{table:MobSuz} in Section~2.4. Specifically, this table gives the values of $\mu_G(H)$ and $|N_G(H)|$ for a set $\mathcal T$ of representatives $H$ of the conjugacy classes of subgroups on which $\mu_G$ can take non-zero values. This information is sufficient for applications of equation~(\ref{NGammaG}): since $|{\rm Aut}\,G|=e|G|$, this now takes the form
\begin{equation}\label{NewNGammaG}
n_{\Gamma}(G)
=\frac{1}{e}\sum_{H\in{\mathcal T}}
\frac{\mu_G(H)|{\rm Hom}(\Gamma,H)|}{|N_G(H)|}.
\end{equation}
For example, it follows that for the Suzuki groups $G$ we have
\[d_2(G)=\frac{1}{e}\sum_{f|e}\mu\left(\frac{e}{f}\right)2^f(2^{4f}-2^{3f}-9),\]
where $\mu$ is the classical M\"obius function on $\mathbb N$, given by $\mu(n)=(-1)^k$ if $n$ is a product of $k$ distinct primes, and $\mu(n)=0$ otherwise. See~\cite{DJ13} for this result, and for a number of other applications of Hall's theory to these groups $G$, and see Section~\ref{Hecke} for an application to certain Hecke groups $\Gamma$.

This paper, together with its applications in~\cite{DJ13}, extends work by Silver and the second author~\cite{JS}, where orientably regular maps of type $\{4,5\}$ with automorphism group $G\cong Sz(q)$ were enumerated, and by Hubard and Leemans~\cite{HL}, where regular maps and polytopes were enumerated. In each case the authors used a restricted form of M\"obius inversion, concentrating mainly on subgroups $H\cong Sz(2^f)$ where $f$ divides $e$. Here we compute $\mu_G(H)$ for all subgroups $H\le G$, allowing the enumeration of a wider range of regular objects, including dessins d'enfants, regular or orientably regular maps and hypermaps,  and coverings of topological spaces (see~\cite{DJ13}).

The method used for calculating the values of $\mu_G$ is as follows. Hall~\cite[Theorem~2.3]{Hal} showed that, in any finite group $G$, a subgroup $H$ satisfies $\mu_G(H)=0$ unless $H$ is an intersection of maximal subgroups of $G$. In our case, with $G=Sz(q)$, instead of directly determining the set $\mathcal I$ of such intersections, we first describe, in Section~\ref{important}, a more convenient set $\mathcal S$ of subgroups of $G$ such that every subgroup in $\mathcal I$ is conjugate to a subgroup in $\mathcal S$ (see Theorem~\ref{conjsgpS}). Since $\mu_G$ is invariant under conjugation, it sufficient to find its values on $\mathcal S$; then the set $\mathcal T$ appearing in equation~(\ref{NewNGammaG}) is simply the subset of $\mathcal S$ on which $\mu_G$ can take non-zero values.

For each pair of subgroups $H, K\in{\mathcal S}$, we determine in Table~2 (Section~5) the number $N(H;K)$ of conjugates of $K$ containing $H$. Since $\mu_G(K)=0$ for all $K\not\in{\mathcal S}$, equation~(\ref{mu}) gives
\begin{equation}\label{recursive}
\mu_G(H)=-\negthickspace\sum_{H<K\in{\mathcal S}}N(H;K)\mu_G(K)
\end{equation}
where $H\in{\mathcal S}\setminus\{G\}$. This allows $\mu_G(H)$ to be calculated recursively, starting with $\mu_G(G)=1$ and then using the values of $\mu_G(K)$ for the subgroups $K\in{\mathcal S}$ properly containing each $H\in{\mathcal S}\setminus\{G\}$.

\medskip

\noindent{\bf Acknowledgements} The authors are grateful to Dimitri Leemans for some very helpful comments on enumeration with Suzuki groups, and to Nikos Kanakis for help in preparing the TeX file.

\section{The Suzuki groups and their subgroups}

This section is based mainly on Suzuki's description of the groups $Sz(q)$ in~\cite{Suz}; see also~\cite[\S XI.3]{HB} and~\cite[\S 4.2]{Wil} for further information. We have largely followed Suzuki's notation for elements and subgroups, except that we use the symbol $F$ for the subgroup denoted in~\cite{Suz} by $H$ (a Frobenius group of order $q^2(q-1)$), while we use $H$ for an arbitrary subgroup of $G$. Also, our rule for distinguishing the subgroups $A_1$ and $A_2$ in \S\ref{defn} is non-standard.

\subsection{The definition of the Suzuki groups}\label{defn}
\label{1.1}
Let ${\mathbb F}={\mathbb F}(e)$ be the finite field ${\mathbb F}_q$ of $q=2^e$ elements for some odd $e \ge 1$, and let $\theta$ be the automorphism $\alpha\mapsto\alpha^r$ of $\mathbb F$ where $r=\sqrt{2q}=2^{(e+1)/2}$, so that $\theta^2$ is the Frobenius automorphism $\alpha \mapsto \alpha^2$.

For any $\alpha, \beta \in \mathbb F$ let $(\alpha, \beta)$ denote the $4 \times 4$ matrix
\[
(\alpha,\beta)=
\left(\,\begin{matrix}
1&&&\cr
\alpha&1&&\cr
\alpha^{\theta+1}+\beta&\alpha^{\theta}&\qquad 1&\cr
\alpha^{\theta+2}+\alpha \beta+\beta^{\theta}&\beta&\qquad\alpha&\qquad 1\cr\end{matrix}\,\right).
\]
Since
$(\alpha, \beta)(\gamma, \delta) = (\alpha + \gamma, \alpha \gamma^{\theta} + \beta + \delta),$
these matrices $(\alpha, \beta)$ form a group $Q=Q(e)$ of order $q^2$, with identity element $(0,0)$.

For each $\kappa\in{\mathbb F}^*:={\mathbb F}\setminus\{0\}$ let $a_{\kappa}$ denote the $4\times 4$ diagonal matrix with diagonal entries $\zeta_i$ where $\zeta_1^{\theta}=\kappa^{1+\theta}$, $\zeta_2^{\theta}=\kappa$, $\zeta_3=\zeta_2^{-1}$ and $\zeta_4=\zeta_1^{-1}$. These matrices form a cyclic group $A_0=A_0(e)\cong {\mathbb F}^*$ of order $q-1$. Since
\[a_{\kappa}^{-1}(\alpha,\beta)a_{\kappa}=(\alpha\kappa,\beta\kappa^{1+\theta}),\]
the group $F=F(e)$ generated by $Q$ and $A_0$ in $GL_4(q)$ is a semidirect product of a normal subgroup $Q$ by a complement $A_0$, so it has order $q^2(q-1)$.

We define $G=G(e)$ to be the subgroup of $GL_4(q)$ generated by $F$ and the $4 \times 4$ matrix with entries $1$ on the minor diagonal and $0$ elsewhere (denoted by $\tau$ in~\cite{Suz}). This is the Suzuki group associated with ${\mathbb F}_q$, usually denoted by $Sz(q)$ or ${}^2B_2(q)$. It is, in fact, the subgroup of the symplectic group $Sp_4(q)=B_2(q)$ fixed by a certain automorphism of order $2$.

In its natural action $g:[v]\mapsto [vg]$ on the projective space ${\mathbb P}^3({\mathbb F})$, $G$  acts as a doubly transitive permutation group of degree $q^2+1$ on the ovoid
\[\Omega = \Omega(e) = 
\{[\alpha^{\theta+2}+\alpha \beta+\beta^{\theta},\,
\beta,\, \alpha,\, 1]\mid \alpha, \beta\in{\mathbb F}\}
\cup \{\infty\} \subset {\mathbb P}^3({\mathbb F}),\]
where
\[\infty:=[1,0,0,0]\in{\mathbb P}^3({\mathbb F}).\]
The subgroup $G_{\infty}$ of $G$ fixing $\infty$ is $F$. This acts as a Frobenius group on $\Omega\setminus\{\infty\}$: its Frobenius kernel is $Q$, acting regularly on $\Omega\setminus\{\infty\}$, and $A_0$ is a Frobenius complement {$G_{\infty,\omega}$, fixing a second point
\[\omega:=[0,0,0,1]\in\Omega\]
and acting semiregularly on $\Omega\setminus\{\infty,\omega\}$. Thus the stabiliser of any three points in $\Omega$ is the identity subgroup $I$, so $G$ acts on $\Omega$ as a Zassenhaus group. 

There are cyclic subgroups of $G$ of mutually coprime odd orders
\[2^e\pm 2^{(e+1)/2}+1=q\pm r+1,\]
contained in Singer subgroups of $GL_4(q)$: note that since $r=\sqrt{2q}$, 
\[(q+r+1)(q-r+1)=q^2+1,\]
which divides $q^4-1$. Let us choose a pair of subgroups $A_1=A_1(e), A_2=A_2(e)$ of $G$ of these two orders, indexed according to the rule
\[|A_1(e)|=a_1(e):=2^e+\chi(e)2^{(e+1)/2}+1,\]
\[|A_2(e)|=a_2(e):=2^e-\chi(e)2^{(e+1)/2}+1,\]
where $\chi(e)=1$ or $-1$ as $e\equiv\pm 1$ or $\pm 3$ mod~$(8)$. (This rule for indexing $A_1$ and $A_2$ differs from that used in~\cite{HB, Suz}, where the rule is that $|A_1(e)|>|A_2(e)|$ for all $e$; the rule adopted here has an advantage which will be explained in \S\ref{subfield}.)

\subsection{Basic properties of Suzuki groups}\label{basic}

Here we record some basic properties of $G$; see~\cite{HB, Suz} for proofs.

\smallskip

\begin{enumerate}

\item $G$ has order $q^2(q^2+1)(q-1)$, and is simple if $e > 1$. (The group $G(1)$ is isomorphic to $AGL_1(5)$, of order $20$.)

\item ${\rm Aut}\,G$ is a semidirect product of ${\rm Inn}\,G\cong G$ by a cyclic group of order~$e$ acting as the Galois group ${\rm Gal}\,{\mathbb F}$ on matrix entries, so $|{\rm Aut}\,G|=e|G|$.

\item Any two subgroups of $G$ conjugate to $Q$ intersect trivially, and any two subgroups conjugate to $F$ have their intersection conjugate to $A_0$.

\item $Q$ is a Sylow $2$-subgroup of $G$ of order $q^2$ and of exponent $4$. The centre $Z$ of $Q$ consists of the identity and the involutions of $Q$ (the matrices $(0,\beta)$ for $\beta\in{\mathbb F}$), with $Z$ and $Q/Z$ elementary abelian of order $q$.

\item $ZA_0\cong F/Z\cong AGL_1(q)$, with $A_0$ acting regularly by conjugation on the non-identity elements of $Z$ and of $Q/Z$.

\item The involutions of $G$ are all conjugate, as are the cyclic subgroups of order $4$; however an element of order $4$ is not conjugate to its inverse.

\item All elements of $G$ except those in a conjugate of $Q$ have odd order. Each maximal cyclic subgroup of $G$ of odd order is conjugate to $A_0$, $A_1$ or $A_2$; the intersection of any two of them is trivial.

\item A non-identity element of $G$ has two fixed points on $\Omega$, one fixed point, or none as it is conjugate to an element of $A_0$, of $Q$ or of $A_i$ for $i=1, 2$, or, equivalently, as it has order dividing $q-1$, $4$ or $q^2+1$.

\end{enumerate}

\subsection{Some particular subgroups}

Here we list some particular subgroups of $G$, in the anticipation that any subgroup $H$ not conjugate to a member of the list satisfies $\mu_G(H)=0$, and can therefore be ignored in the enumerations mentioned in Section~\ref{Intro}. This list enables us to give a statement of all values of the M\"obius function, prior to proving it.

\subsubsection{Subgroups associated with subfields}\label{subfield}

If $f$ divides $e$ then restricting matrix entries to the subfield ${\mathbb F}(f)$ of $\mathbb F$ of order $2^f$ yields a subgroup $G(f) = Sz(2^f)$ of $G$. This acts doubly transitively, with degree $2^{2f}+1$, on the subset $\Omega(f)$ of $\Omega$ defined over ${\mathbb F}(f)$. Since the point $\infty$ is defined over the prime field ${\mathbb F}(1)$, its stabiliser in $G(f)$ is $F(f):=F\cap G(f)$, which acts faithfully on $\Omega(f)\setminus\{\infty\}$ as a Frobenius group with kernel $Q(f):=Q\cap G(f)=\{(\alpha,\beta)\mid\alpha,\beta\in{\mathbb F}(f)\}$ and complement $A_0(f):=A_0\cap G(f)=\{a_{\kappa}\mid\kappa\in{\mathbb F}(f)^*\}$. (Note that since $r$ and $r+1$ are coprime to $q-1$, we have $a_{\kappa}\in G(f)$ if and only if $\kappa\in{\mathbb F}(f)^*$.) Let $Z(f)$ denote the centre of $Q(f)$, an elementary abelian group of order $2^f$. If $f$ and $h$ are divisors of $e$, and $f$ divides $h$, then
\[G(f)\leq G(h), \; F(f)\leq F(h), \; Q(f)\leq Q(h), \; Z(f)\leq Z(h),\; A_0(f)\leq A_0(h).\]
Indeed, each of these five families of subgroups forms a lattice isomorphic to the lattice $\Lambda(e)$ of divisors of $e$, a fact which is useful in evaluating $\mu_G$.

\medskip

Similarly, we would like to choose cyclic subgroups $A_i(f)$ of $G$ of order
\[|A_i(f)|=a_i(f)=2^f\pm\chi(f)2^{(f+1)/2}+1\]
corresponding to $A_i$ for each $i=1, 2$ so that they also satisfy $A_i(f)\le A_i(h)$ whenever $f$ divides $h$. The following result allows this:

\begin{lemma}\label{indexing}
For each $i=1, 2$ and each $f$ dividing $e$ there is a cyclic subgroup $A_i(f)$ of $A_i$ of order $a_i(f)$.
\end{lemma}

\noindent{\sl Proof.} Since $A_i$ is cyclic, it is sufficient to show that if $f$ divides $e$ then $a_i(f)$ divides $a_i(e)$ for $i=1, 2$. Let $m:=2^{2f}+1=a_1(f)a_2(f)$. Then $2^j\equiv -2^{j-2f}$ mod~$(m)$ for all $j\ge 2f$, so
\[2^e\equiv -2^{e-2f}\equiv 2^{e-4f}\equiv\cdots\equiv(-1)^{(e-f)/2f}2^f\;{\rm mod}~(m)\]
since $e$ is odd. Now $(e-f)/2f$ is even or odd as $e/f\equiv \pm 1$ mod~$(4)$, giving $2^e\equiv\pm 2^f$ mod~$(m)$ respectively. Similarly,
\[2^{(e+1)/2}\equiv (-1)^k\,2^l\;{\rm mod}~(m)\]
where
\[k=\lfloor\frac{e+1}{4f}\rfloor
\quad{\rm and}\quad
0\le l=\frac{e+1}{2}-2kf<2f.\]
Now $k$ is even if $e/f\equiv 1$ or $3$ mod~$(8)$, and odd if $e/f\equiv -1$ or $-3$ mod~$(8)$, while $l=(f+1)/2$ or $(3f+1)/2$ as $e/f\equiv \pm 1$ mod~$(4)$.

In the case where $e/f\equiv 1$ mod~$(8)$, so that $\chi(e)=\chi(f)$, it follows that
\[a_i(e)=2^e\pm\chi(e)2^{(e+1)/2}+1\equiv 2^f\pm\chi(f)2^{(f+1)/2}+1=a_i(f)\;{\rm mod}~(m)\]
for $i=1, 2$; since $a_i(f)$ divides $m$ it divides $a_i(e)$. The argument is essentially the same if $e/f\equiv -3$ mod~$(8)$, with a change of sign in the coefficient of $2^{(f+1)/2}$ balanced by a change of sign in $\chi(f)$. If $e/f\equiv 3$ or $-1$ mod~$(8)$ a similar argument, using the factorisation
\[-2^f\pm 2^{(3f+1)/2}+1=(1\pm 2^{(f+1)/2})(2^f\mp 2^{(f+1)/2}+1),\]
gives the result.\hfill$\square$

\medskip

In fact, since $A_i$ is cyclic, this shows that if $f \mid h \mid e$ then $A_i(f)\le A_i(h)$, so for each $i=1, 2$ the subgroups $A_i(f)$, where $f$ divides $e$, form a lattice isomorphic to $\Lambda(e)$. This explains our non-standard choice of indexing for these two sets of subgroups. Note, however, that $A_i(f)$ is now not necessarily a subgroup of $G(f)$, though it is conjugate to such a subgroup.

\subsubsection{The normalisers of some subgroups}\label{normalisers}

The normaliser $B_0$ of $A_0$ in $G$ is a dihedral group of order $2(q-1)$; it is the subgroup $G_{\{\infty,\omega\}}$ of $G$ preserving the subset $\{\infty, \omega\}$ of $\Omega$, with its subgroup $A_0$ fixing these two elements and its involutions transposing them. Let us choose a particular involution $c\in B_0$ and, for each $f$ dividing $e$, define \[B_0(f):=\langle A_0(f),c\rangle\le B_0,\]
a dihedral group of order $2(2^f-1)$ (so $B_0(1)\cong C_2$). If $f > 1$, then $B_0(f)$ is 
self-normalising whereas the normaliser of $A_0(f)$ is $B_0$.

For $i=1, 2$ the normaliser $B_i$ of $A_i$ in $G$ is a semidirect product of $A_i$ and a cyclic group of order $4$ generated by an element $c_i$ satisfying $c_i^{-1} a c_i = a^{2^{e}}$ for all $a \in A_i$.
For each $f$ dividing $e$ let
\[B_i(f):=\langle A_i(f), c_i\rangle\le B_i,\]
so $|B_i(f)| = 4a_i(f)$, with $B_2(1)\cong C_4$. If $i=1$ or $f>1$ then $B_i(f)$ is self-normalising, whereas the normaliser of $A_i(f)$ is $B_i$. 

By their construction, these groups $B_i(f)$ are (abstract) Frobenius groups of degree $a_i(f)$, and they satisfy $B_i(f)\le B_i(h)$ for $i=0, 1$ and $2$ whenever $f\mid h\mid e$.

\bigskip

\subsubsection{An important set of subgroups}\label{important}

For each $f$ dividing $e$, we have defined the following subgroups of $G$, with the symbols $(f)$ usually omitted when $f=e$:
\begin{equation}\label{setS}
G(f),\; F(f),\; Q(f),\; Z(f),\; B_i(f),\; A_i(f)\quad (i=0, 1, 2).
\end{equation}
Let $\mathcal S$ denote the set consisting of the subgroups in~(\ref{setS}) for all $f$ dividing $e$. The conjugacy class in $G$ of any of these groups will be denoted by changing the appropriate italic capital letter to the corresponding script capital; thus ${\mathcal G}(f), {\mathcal F}, \ldots$ denote the conjugacy classes containing $G(f)$, $F$, and so on.

We note the following coincidences{, conjugacies (denoted by $\sim$) and isomorphisms:
\[G(1)\sim B_1(1) \cong AGL_1(5),\quad F(1)=Q(1)\sim B_2(1)\cong C_4,\]
\[B_0(1)\sim Z(1)\cong C_2, \quad A_2(1)=A_0(1)=I.\]
In addition, if $3$ divides $e$ then
\[B_1(1)=B_1(3),\quad A_1(1)=A_1(3)\cong C_5.\]
Apart from these, any two distinct terms in~(\ref{setS}) represent non-conjugate subgroups of $G$. In view of their special role in the following calculations, we will denote the class ${\mathcal A}_2(1)={\mathcal A}_0(1)$ by ${\mathcal C}_1$,  the class ${\mathcal B}_0(1)={\mathcal Z}(1)$ by ${\mathcal C}_2$, and the class ${\mathcal F}(1)={\mathcal Q}(1)={\mathcal B}_2(1)$ by ${\mathcal C}_4$, since these consist of the cyclic subgroups of $G$ of orders $1, 2$ and $4$.

\bigskip

\subsection{The M\"obius function of a Suzuki group}\label{MobSuz}

We can now present the main result of this paper in the form of Table~\ref{table:MobSuz}, which gives the non-zero values of $\mu_G(H)$ for the subgroups $H$ of $G$; any subgroups $H$ not appearing in Table~\ref{table:MobSuz} (such as $Q(f)$ and $Z(f)$ for $f>1$) satisfy $\mu_G(H)=0$, and can therefore be ignored in applying equations such as~(\ref{NGammaG}). Because of the conjugacies listed in \S2.3.3, some conjugacy classes appear `under an alias': for instance ${\mathcal F}(1)$ appears as ${\mathcal B}_2(1)$, and if $3$ divides $e$ then ${\mathcal G}(1)$ and ${\mathcal B}_1(1)$ appear as ${\mathcal B}_1(3)$. In the second column, $a_i(f)=2^f\pm\chi(f)2^{(f+1)/2}+1$ for $i=1,2$ (see \S\ref{defn} and \S\ref{subfield}). In the third column, the values of $|N_G(H)|$ are given for applications of equation~(\ref{NewNGammaG}). In the final column, $\mu$ is the classical M\"obius function on $\mathbb N$, defined by
\[\sum_{m|n}\mu(f)=\delta_{n,1}\]
for all $n\in{\mathbb N}$, with the consequence that $\mu(n)=(-1)^k$ or $0$
as $n$ is or is not a product of $k$ distinct primes for some $k\ge 0$.

\begin{table}[ht]
\centering
\begin{tabular}{| p{2.9cm} | p{3.5cm} | p{1.6cm} | p{2.4cm} |}
\hline
Conjugacy class of $H$ & $|H|$ &  $|N_G(H)|$ & $\mu_G(H)$ \\
\hline\hline
$\mathcal{G}(f),\; \; \, 1<f\mid e$ & $2^{2f}(2^{2f}+1)(2^f-1)$ & $|H|$ & $\mu(e/f)$\\
\hline
$\mathcal{F}(f),\; \; 1<f\mid e$ & $2^{2f}(2^f-1)$ & $|H|$ & $-\mu(e/f)$ \\
\hline
$\mathcal{B}_0(f),\; 1<f\mid e$ & $2(2^f-1)$ & $|H|$ & $-\mu(e/f)$ \\
\hline
$\mathcal{A}_0(f),\, 1<f\mid e$ & $2^f-1$ & $2(q-1)$ & $2\frac{(2^e-1)}{(2^f-1)}\mu(e / f)$ \\ [0.6ex]
\hline
$\mathcal{B}_1(f),\; 1<f\mid e$ & $4a_1(f)$ & $|H|$ & $-\mu(e/f)$ \\
\hline
$\mathcal{B}_2(f),\; 1<f\mid e$ & $4a_2(f)$ & $|H|$ & $-\mu(e/f)$ \\
\hline
${\mathcal B}_2(1)=\mathcal{C}_4$ & $4$ & $2q$ & $-2^{e}\mu(e)$ \\
\hline
${\mathcal B}_0(1)=\mathcal{C}_2$ & $2$ & $q^2$ & $-2^{2e-1}\mu(e)$ \\
\hline
${\mathcal A}_0(1)=\mathcal{C}_1$ & $1$ & $|G|$ & $|G|\mu(e)$ \\
\hline

\end{tabular}
\caption{Information about the subgroups $H$ with non-zero values for $\mu_G(H)$}
\label{table:MobSuz}
\end{table}

Our aim is to show that the final column of this table is correct, by proving the following theorem:

\begin{thm}\label{mainthm}
Let $G$ be a Suzuki group $Sz(2^e)$ for some odd $e>1$, and let $H$ be a subgroup of $G$. If $\mu_G(H)\ne 0$ then $\mu_G(H)$ is as given by Table~1.
\end{thm}

Note that when $H$ is the identity subgroup, the value of $\mu_G(H)$ is $|G|\mu(e)$, supporting a conjecture of Conder that this value is divisible by $|G|$ whenever $G$ is an almost simple group, that is, $S\le G\le{\rm Aut}\,S$ for some non-abelian finite simple group $S$. By~\cite{DowJLMS, Hal} this is true when $G=L_2(q)$ or $PGL_2(q)$, and computer calculations have verified the conjecture when $S$ is one of the smaller alternating or sporadic simple groups.

\section{Subgroups $H$ with $\mu_G(H)\ne 0$}

As a first step towards proving Theorem~\ref{mainthm}, in this section we find some necessary conditions for a subgroup $H$ of $G$ to satisfy $\mu_G(H)\ne 0$.

\subsection{Maxint subgroups}\label{maxint}

If $G$ is any finite group, we shall say that a subgroup $H$ of $G$ is {\sl maxint\/} if it is the intersection of a set of maximal subgroups of $G$ (when $H=G$ this set is empty). The set of maxint subgroups of $G$ will be denoted by $\mathcal I$. Hall proved in~\cite[Theorem~2.3]{Hal} that if $H \notin{\mathcal I}$ then $\mu_G(H) = 0$, so in determining $\mu_G$ one may restrict attention to the subgroups $H\in{\mathcal I}$.  Since $\mu_G$ is preserved under conjugacy, it is sufficient to consider a set of representatives of the conjugacy classes of subgroups in $\mathcal I$. The main step in the proof of Theorem~\ref{mainthm} is to show that if $G$ is a Suzuki group $G(e)$ then the set $\mathcal S$ defined in \S\ref{important} contains such a set of representatives:

\begin{thm}\label{conjsgpS}
If $H \in{\mathcal I}$ then $H$ is conjugate in $G$ to a subgroup in $\mathcal S$, that is, $\mathcal I$ is contained in the union of the conjugacy classes
\[{\mathcal G}(f),\;{\mathcal F}(f),\;{\mathcal Q}(f),\;{\mathcal Z}(f),
\;{\mathcal B}_i(f),\;{\mathcal A}_i(f)\]
 of subgroups of $G$, where $f$ divides $e$ and $i=0, 1$ or $2$.
\end{thm}

The rest of this section is devoted to a proof of this theorem. We will use the following criterion for a subgroup $H$ of $G$ to be in $\mathcal I$. Let $\mathcal M$ denote the set of all maximal subgroups of $G$, and let ${\mathcal M}(H)$ denote the set of those containing a particular subgroup $H$ of $G$. Then the following, valid for any finite group $G$, is evident:

\begin{lemma}\label{Itest}
Let $H\le G$. Then
\[H\;\le\negthickspace\bigcap_{M\in{\mathcal M}(H)}\negthickspace M\,,\]
with equality if and only if $H\in{\mathcal I}$.
\end{lemma}

\subsection{Maximal subgroups}\label{maxsgps}

We will systematically apply Lemma~\ref{Itest} to the various subgroups $H$ of $G$, using the following result:

\begin{prop}\label{max}
The set $\mathcal M$ of maximal subgroups of $G$ is given by
\[{\mathcal M}\;=\negthickspace\bigcup_{e/f\,{\rm prime}}\negthickspace\negthickspace{\mathcal G}(f)\cup{\mathcal F}\cup{\mathcal B}_0\cup{\mathcal B}_1\cup{\mathcal B}_2.\]
\end{prop}

This result is an immediate consequence of the following classification, due to Suzuki~\cite[Theorems 9 and 10]{Suz}:

\begin{prop}\label{subgps}
If $H\le G$ then either $H\in{\mathcal G}(f)$ for some $f$ dividing $e$, or $H$ is a subgroup of a group in $\mathcal F$ or in ${\mathcal B}_i$ for some $i=0, 1$ or $2$.
\end{prop}

In the first case $H$ is either simple or isomorphic to $G(1)\cong AGL_1(5)$, and in the second case $H$ is solvable. Finite solvable groups $H$ all satisfy Hall's theorems~\cite{Hal28} on the existence and conjugacy of Hall $\pi$-subgroups for any set $\pi$ of primes, generalising Sylow's theorems for single primes. We will use this fact, mainly with $\pi$ the set $2'$ of odd primes.

Since $G(f)\in{\mathcal S}$ for each $f$ dividing $e$, it follows from Proposition~\ref{subgps} that, in proving Theorem~\ref{conjsgpS}, it is sufficient to assume that $H$ is a subgroup of a group in $\mathcal F$ or ${\mathcal B}_i$ for $i=0, 1$ or $2$. We will deal with these cases in turn.

In preparation for applying Lemma~\ref{Itest} in the first case, we will consider how the various maximal subgroups of $G$ intersect $F$.

\subsection{Point-stabilisers in maximal subgroups}
\label{intersections}

Recall that $F$ is the stabiliser in $G$ of the point $\infty\in\Omega$. If $H\le F$ then
\[\negthickspace\bigcap_{M\in{\mathcal M}(H)}\negthickspace M
=\negthickspace\bigcap_{M\in{\mathcal M}(H)}\negthickspace (M\cap F),\]
so in applying Lemma~\ref{Itest} to $H$ one can restrict attention to the stabilisers $M_{\infty}=M\cap F$ of $\infty$ for the various maximal subgroups $M$ of $G$. The following result describes the possibilities for these stabilisers.

\begin{lemma}\label{McapF}
Let $M$ be a maximal subgroup of $G$. 
\begin{enumerate}
\item If $M=F^g\in{\mathcal F}$, then $M\cap F=F$ or $M\cap F=G_{\infty,\infty g}\in{\mathcal A}_0$ as $g\in F$ or not.
\item If $M=G(f)^g\in{\mathcal G}(f)$ for some $f|e$, then $M\cap F=F(f)^g\in{\mathcal F}(f)$ or $M\cap F=I\in {\mathcal C}_1$ as $g\in G(f)F$ or not.
\item If $M\in{\mathcal B}_0$ then $M\cap F\in{\mathcal A}_0\cup{\mathcal C}_1\cup{\mathcal C}_2$.
\item If $M\in{\mathcal B}_i$ for $i=1,2$ then $M\cap F\in{\mathcal C}_1\cup{\mathcal C}_2\cup{\mathcal C}_4$.
\end{enumerate}
\end{lemma}

In order to prove part~(2), we first need the following lemma:

\begin{lemma}\label{semireg}
If $f$ divides $e$ then $G(f)$ acts semi-regularly on $\Omega\setminus\Omega(f)$.
\end{lemma}

\noindent{\sl Proof.} By \S\ref{basic}(8), a non-identity element of $G$ fixes $2, 1$ or $0$ elements of $\Omega$ as it has order dividing $q-1$, $4$ or $q^2+1$. Similarly, a non-identity element of $G(f)$ fixes $2, 1$ or $0$ elements of $\Omega(f)$ as it has order dividing $2^f-1$, $4$ or $2^{2f}+1$ respectively. Since $2^f-1$ divides $q-1$, and $2^{2f}+1$ divides $q^2+1$, a non-identity element of $G(f)$ can have no further fixed points in $\Omega\setminus\Omega(f)$. Thus all orbits of $G(f)$ on this set are regular, with point-stabilisers $G(f)\cap G_{\alpha}=I$ for $\alpha\in\Omega\setminus\Omega(f)$.\hfill$\square$

\bigskip

\noindent{\sl Proof of Lemma~\ref{McapF}.} The maximal subgroups $M$ of $G$ are given by Proposition~\ref{max}.

(1) This part is trivial, since $F$ and $M$ are the stabilisers in $G$ of $\infty$ and $\infty g$, and $G$ is doubly transitive on $\Omega$.

(2) If $M=G(f)^g\in{\mathcal G}(f)$ then Lemma~\ref{semireg} shows that $M$ acts doubly transitively on $\Omega(f)g$, and semiregularly on its complement. Thus $M\cap F=F(f)^g$ or $I$ as $\infty\in\Omega(f)g$ or not, that is, as $g\in G(f)F$ or not.

(3) Each $M\in{\mathcal B}_0$ is the subgroup $G_{\{\alpha,\beta\}}$ of $G$ preserving an unordered pair $\{\alpha, \beta\}\subset\Omega$. If $\infty\not\in\{\alpha, \beta\}$ then since $G_{\alpha,\beta,\infty}=I$ we have $|M\cap F|\le 2$, whereas if $\infty=\alpha$ or $\beta$ then $M\cap F=G_{\alpha,\beta}\in{\mathcal A}_0$.

(4) If $M\in{\mathcal B}_i$ for $i=1$ or $2$ then $M=N_G(A)\cong A\rtimes C_4$ for some $A\in{\mathcal A}_i$; since $A$ acts without fixed points, $M\cap F$ is isomorphic to a subgroup of $C_4$, so it is in ${\mathcal C}_m$ for $m=1, 2$ or $4$.\hfill$\square$

\subsection{Subgroups $H$ of $F$}

We first prove Theorem~\ref{conjsgpS} for subgroups $H\in{\mathcal I}$ which are contained in groups in $\mathcal F$. Replacing $H$ with a conjugate, we may assume that $H\le F$.

\subsubsection{Preliminaries}\label{prelim}

Here we record some observations and simplifications which will be used in the proof.

\smallskip

(a) Lemma~\ref{McapF} shows that each $M\in{\mathcal M}(H)$ satisfies $M\cap F\in{\mathcal X}_M$ where ${\mathcal X}_M={\mathcal F}(f_M)$ for some $f_M$ dividing $e$ (depending on $M$), or ${\mathcal X}_M={\mathcal A}_0$, or ${\mathcal X}_M={\mathcal C}_m$ for some $m$ dividing $4$. If ${\mathcal X}_M={\mathcal C}_m$ for some $m$, then since $H\le M\cap F$ we have $H\in{\mathcal C}_1\cup{\mathcal C}_2\cup{\mathcal C}_4$, so $H$ is as required, i.e.~conjugate to an element of $\mathcal S$; we may therefore assume that for each $M\in{\mathcal M}(H)\setminus\{F\}$ we have ${\mathcal X}_M={\mathcal F}(f_M)$ for some $f_M$ or ${\mathcal X}_M={\mathcal A}_0$, with $M\in{\mathcal G}(f_M)$ or $\mathcal F$ respectively.

\smallskip

(b) Since $F$ is solvable, Hall's theorems~\cite{Hal28} imply that a Hall $2'$-subgroup $A$ of $H$ is contained in one of the Hall $2'$-subgroups of $F$. These are the conjugates of $A_0$, and $Q$  permutes them regularly by conjugation, so by conjugating $H$ with a suitable element of $Q$ we may assume that $A\le A_0$.

\smallskip

(c) If $H$ has even order it contains an involution. The involutions in $F$ (the non-identity elements of $Z$) are all conjugate under $A_0$, so in this case, by conjugating $H$ with an element of $A_0$ we may also assume that $H$ contains $z:=(0,1)$.

\smallskip

(d) If any $M\in{\mathcal M}(H)\setminus\{F\}$ satisfies $M=G(f_M)^g\in{\mathcal G}(f_M)$ for some {$g\in G$, then since $M\cap F\in{\mathcal F}(f_M)$ we have $g\in G(f_M)F$ by Lemma~\ref{McapF}. Without loss of generality we can therefore choose this conjugating element $g$ to be in $F$. Then
\[F(f_M)^g=(G(f_M)\cap F)^g=G(f_M)^g\cap F=M\cap F.\]
Thus $z\in M\cap F$, so $z^{g^{-1}}\in F(f_M)$. Since $F=A_0Q$ we can write $g=ab$ where $a\in A_0$ and $b\in Q$. Since $z$ is in the centre $Z$ of $Q$ we have
\[z^{g^{-1}}=(z^{b^{-1}})^{a^{-1}}=z^{a^{-1}},\]
so $z^{a^{-1}}\in F(f_M)$. Since $z\in F(1)$ and $A_0$ acts regularly by conjugation on the involutions in $Z$, this implies that $a\in A_0(f_M)$. Thus $g=ab$ with $a\in F(f_M)$, so each $M\in{\mathcal M}(H)\setminus\{F\}$ satisfies
\[M\cap F=F(f_M)^g=F(f_M)^b\]
for some $f_M$ dividing $e$, with $b\in Q$.

\smallskip

(e) We claim that if $f\mid h\mid e$ then the set
\[Q(f,h):=\{g\in Q \mid Q(f)^g\le Q(h)\}\]
is the union of the cosets $(\alpha,0)Z$ of $Z$ in $Q$ where $\alpha\in{\mathbb F}(h)$.

Clearly this set consists of complete cosets of $Z$ in $Q$. The elements $(\alpha,0)$ where $\alpha\in{\mathbb F}$ are representatives of these cosets, since there is an epimorphism $(\alpha,\beta)\mapsto\alpha$ from $Q$ to the additive group of $\mathbb F$, with kernel $Z$. Therefore it suffices to show that $g:=(\alpha,0)\in Q(f,h)$ if and only if $\alpha\in{\mathbb F}(h)$.

If $\alpha\in{\mathbb F}(h)$ then $g\in Q(h)$; since $Q(f)\le Q(h)$ we have $Q(f)^g\le Q(h)$ and hence $g\in Q(f,h)$. For the converse, note that $(1,0)\in Q(f)$. A simple calculation shows that
\[(1,0)^g=(1,\alpha+\alpha^{\theta}),\]
so if $g\in Q(f,h)$ then $\alpha+\alpha^{\theta}\in{\mathbb F}(h)$. The function  $\phi:x\mapsto x+x^{\theta}$ maps each subfield $\mathbb K$ of ${\mathbb F}$ into itself. Composing $\phi$ with itself gives a quadratic polynomial 
\[\phi^2:x\mapsto (x+x^{\theta})+(x+x^{\theta})^{\theta}=x+x^{\theta^2}=x+x^2\]
defined over the prime field, so if $\beta\in{\mathbb K}$ then any element of $\phi^{-2}(\beta)$ has degree at most $2$ over $\mathbb K$. Since $e$ is odd, $\mathbb F$ contains no quadratic extensions, so $\phi^{-2}({\mathbb K})\subseteq{\mathbb K}$ and hence $\phi^{-1}({\mathbb K})\subseteq{\mathbb K}$. In particular, since $\phi(\alpha)\in{\mathbb F}(h)$ we have $\alpha\in{\mathbb F}(h)$.

\medskip

We can now start the case-by-case analysis of maxint subgroups $H\le F$.

\subsubsection{Subgroups $H\le F$ which are not $2$-groups}

First suppose that $H$ is not a $2$-group, or equivalently the Hall $2'$-subgroup $A\le A_0$ of $H$ is not the identity subgroup, so that $C_G(A)=A_0$. For each $M\in{\mathcal M}$, $A$ is contained in a maximal cyclic subgroup $A_M$ of $M$, which has order $2^{f_M}-1$ if $M\in{\mathcal G}(f)$, and order $2^e-1$ if $M\in{\mathcal F}$. This subgroup $A_M$ centralises $A$, so it is contained in $A_0$; thus $A=A_0(f_M)$, where we take $f_M=e$ when $M\in{\mathcal F}$, since this is the unique subgroup of $A_0$ of order $2^{f_M}-1$. Now $H$ is the intersection of these subgroups $M$, so $A$ is the intersection of the corresponding subgroups $A_M$; it therefore has the form $A=A_0(f)$ where $f$ is the highest common factor of the divisors $f_M$ of $e$.

If $H$ has odd order then this gives $H=A\in{\mathcal A}_0(f)$, one of the types allowed for in Theorem~\ref{conjsgpS}, so we may assume that $H$ has even order. As noted in \S\ref{prelim}(c), this allows us to assume that $z\in H$. This also implies that each $M\in{\mathcal M}(H)\setminus\{F\}$ is in ${\mathcal G}(f_M)$ for some $f_M$, for otherwise $M\in{\mathcal F}$ and so $H$ is a subgroup of a two-point stabiliser $M\cap F$, which has odd order. As shown in \S\ref{prelim}(d), it follows that such subgroups $M$ satisfy $M\cap F=F(f_M)^b$ for some $b\in Q$.

We have $A_0(f_M)=A\le A_M\le M\cap F=F(f_M)^b$, so $A_0(f_M)$ and $A_0(f_M)^{b^{-1}}$ are both point-stabilisers in $F(f_M)$; because $F(f_M)$ acts as a Frobenius group on $\Omega(f)$, its kernel $Q(f_M)$ permutes these point-stabilisers regularly by conjugation, so $A_0(f_M)=A_0(f_M)^{b^{-1}c}$ for some $c\in Q(f_M)$. Thus the element $b^{-1}c$ of $Q$ normalises $A_0(f_M)$, so it also normalises $C_G(A_0(f_M))=A_0$. However, $Q$ permutes the conjugates of $A_0$ regularly by conjugation (since it is also a Frobenius group), so $b^{-1}c=1$ and hence $b=c\in Q(f_M)$.

Thus $M\cap F=F(f_M)^b=F(f_M)$ for each $M\in{\mathcal M}(H)\setminus\{F\}$, so $H$, being the intersection of such subgroups $F(f_M)$, together with $F$, has the form $F(f)\cap F=F(f)$ for some divisor $f$ of $e$, giving $H\in{\mathcal F}(f)$ as required.

\subsubsection{Subgroups $H\le F$ which are $2$-groups}

Now suppose that $H$ is a $2$-group, so $H\le Q$. By \S\ref{prelim}(a), for each $M\in{\mathcal M}(H)\setminus\{F\}$ either $M\cap F\in{\mathcal F}(f_M)$ for some $f_M$ dividing $e$, or $M\cap F\in{\mathcal A}_0$, with $M\in{\mathcal G}(f_M)$ or $\mathcal F$ respectively. We may assume that $H\ne I$, so $H$ has even order and hence (as in \S3.4.2) the second possibility cannot arise. Thus $M\cap F=F(f_M)^b$ for some $b\in Q$, as shown in \S\ref{prelim}(d). As $Q$ is normal in $F$, and is a Sylow $2$-subgroup of $F$, we have $M\cap Q=Q(f_M)^b$; comparing centres, we see that $M\cap Z=Z(f_M)^b=Z(f_M)$ since $Z(f_M)$, being central in $Q$, is normalised by $b$. Thus if $H\le Z$ then
\[H \; =\bigcap_{M\in{\mathcal M}(H)}\negthickspace\negthickspace M \;
=\bigcap_{M\in{\mathcal M}(H)}\negthickspace\negthickspace(M\cap Z) \;
=\bigcap_{M\in{\mathcal M}(H)}\negthickspace\negthickspace Z(f_M) \;
= \; Z(f),\]
where $f$ is the highest common factor of the integers $f_M$, so $H\in{\mathcal Z}(f)$.

We may therefore assume that $H\not\le Z$, so $H$ contains an element of order $4$. Since $F$ has a single conjugacy class of cyclic subgroups of order $4$, we may assume that $H$ contains the subgroup $Q(1)=\{1, z, y^{\pm 1}\}$ where $y:=(1,0)$. Then $Q(1)\le H\le M\cap Q=Q(f_M)^b$, so by \S3.4.1(e) we have $b^{-1}\in (\alpha,0)Z$ for some $\alpha\in{\mathbb F}(f_M)$. This shows that $M\cap Q=Q(f_M)^b=Q(f_M)$ for each $M\in{\mathcal M}(H)\setminus\{F\}$, so taking the intersection over all such $M$ gives $H=Q(f)\in{\mathcal Q}(f)$ where $f$ is the highest common factor of the integers $f_M$.

\subsection{Subgroups $H$ of $B_i$}

Now suppose that $H$ is a subgroup of a group in ${\mathcal B}_i$ for some $i=0, 1$ or $2$, and is maxint. Without loss of generality we may assume that $H\le B_i$.

\subsubsection{Subgroups $H$ of $B_0$}

Suppose that $H\le B_0$. The subgroup $A:=H\cap A_0=H\cap F$ is maxint, since $H$ is, it is contained in $F$, and it has odd order, so by an argument in \S3.4.2 we see that $A=A_0(f)$ for some $f$ dividing $e$. Now $H$ contains $A$ with index at most $2$, so either $H=A_0(f)$, or $H$ is a dihedral subgroup of $B_0$ conjugate (since $|A_0|$ is odd) to $B_0(f)$. Thus $H$ is in ${\mathcal A}_0(f)$ or ${\mathcal B}_0(f)$.

\subsubsection{Subgroups $H$ of $B_i$ for $i=1$ or $2$}

Suppose that $H\le B_i$ where $i=1$ or $2$. Let $A:=H\cap A_i$. If $|A|=1$ then $H$ is isomorphic to a subgroup of $B_i/A_i\cong C_4$, so $H$ is in ${\mathcal C}_m$ for some $m=1, 2$ or $4$. We may therefore assume that $|A|>1$. Any subgroup $M\in{\mathcal M}(H)\setminus\{B_i\}$ contains $A$, which has order dividing $q^2+1$, so it follows from the classification of the maximal subgroups in Proposition~\ref{max} that $M$ must be in ${\mathcal B}_i$ or in ${\mathcal G}(f_M)$ for some $f_M$ dividing $e$. The first possibility can be dismissed, since distinct subgroups in ${\mathcal B}_i$ have intersections of order dividing $4$, so $M\in{\mathcal G}(f_M)$. We can now argue as in \S3.4.2, by considering subgroups centralising $A$, to show that $A=A_i(f)$ for some $f$ dividing $e$. 

Now $|H:A|$ divides $|B_i:A_i|=4$. If $|H:A|=1$ then $H=A_i(f)\in{\mathcal A}_i(f)$, as required. If $|H:A|=4$ then $H$ is a subgroup of $B_i$ of order $4a_i(f)$; all subgroups of this order are conjugate in $B_i$ to $B_i(f)$, so $H\in{\mathcal B}_i(f)$. We will show that the remaining case $|H:A|=2$, where $|H|=2a_i(f)$, cannot arise.

In a Suzuki group $Sz(q)$, any subgroup $K$ of order $2m$, where $m$ divides $q^2+1$, is contained in a unique subgroup $K^*$ of order $4m$. (This is because $A_i<KA_i<B_i$ for $i=1$ or $2$, up to conjugacy, and the complements for $A_i$ in $B_i$ have mutually trivial intersections.) Applying this to the subgroup $K=H$, firstly as a subgroup of $G$, and then as a subgroup of each of the Suzuki subgroups $M\cong G(f)$ in ${\mathcal M}(H)$, we see that there is a subgroup $H^*\le B_i$, containing $H$ with index $2$, such that $H^*\le M$ for all $M\in{\mathcal M}(H)$. Lemma~\ref{Itest} then shows that $H\not\in{\mathcal I}$

This completes the proof of Theorem~\ref{conjsgpS}. \hfill$\square$

\section{Size of conjugacy classes}\label{conjclasses}

An important step in proving the statement of the M\"obius function of $G$ in Theorem~\ref{mainthm} is to determine the number of conjugates of each subgroup $H\in\mathcal S$, equal to the index in $G$ of its normaliser $N_G(H)$. The orders of some of these normalisers are noted in Table~1.

\begin{thm}\label{conjsize}
Let $f$ divide $e$. Then
\begin{enumerate}
\item $N_G(G(f))=G(f)$ and $|\mathcal{G}(f)| = |G|/|G(f)|$;
\item $N_G(F(f))=F(f)$ and $|\mathcal{F}(f)| = |G|/|F(f)|$ if $f>1$;
\item $|N_G(Q(f))|=2^{e+f}(2^f-1)$ and $|\mathcal{Q}(f)| = |G|/2^{e+f}(2^f-1)$ if $f>1$;
\item $N_G(Z(f))=QA_0(f)$ and $|\mathcal{Z}(f)| = |G|/2^{2e}(2^f-1)$ if $f>1$;
\item $N_G(B_i(f))=B_i(f)$ and $|\mathcal{B}_i(f)| = |G|/|B_i(f)|$ if $i=1$, or if $i=0$ or $2$ and $f>1$;
\item $N_G(A_i(f))=B_i$ and $|\mathcal{A}_i(f)| = |G|/|B_i|$ if $i=1$, or if $i=0$ or $2$ and $f>1$;
\item $|N_G(B_2(1))|=2q$ and $|\mathcal{B}_2(1)|=q(q^2+1)(q-1)/2$;
\item $|N_G(B_0(1))|=q^2$ and $|\mathcal{B}_0(1)|=(q^2+1)(q-1)$.
\end{enumerate}
\end{thm}

\noindent{\sl Proof.} (1) Let $H=G(f)$ where $f$ divides $e$. If $f>1$, then since $N_G(G(f))$ contains $G(f)$ it cannot be solvable, so by Proposition~\ref{subgps} it must be conjugate to $G(h)$ for some multiple $h$ of $f$. Since $G(h)$ is simple, we must have $h=f$ and $N_G(H)=H$, giving $|{\mathcal G}(f)|=|G|/|H|$. The case $f=1$ is dealt with in (5), since $G(1)$ is conjugate to $B_1(1)$.

\smallskip

\noindent(4) It is convenient to prove (4) before (2) and (3). Let $f>1$. Any element of $G$ normalising $Z(f)$ must fix its unique fixed point $\infty$, so $N_G(Z(f))\le F$. By \S\ref{basic}(4), $F=QA_0$. Now $Z(f)$ is centralised by $Q$ since it lies in the centre $Z$ of $Q$, and \S\ref{basic}(5) implies that  $N_G(Z(f))\cap A_0=A_0(f)$, so $N_G(Z(f))=QA_0(f)$, of order $|Q|.|A_0(f)|=q^2(2^f-1)=2^{2e}(2^f-1)$.

\smallskip

\noindent(3) Any element of $G$ normalising $Q(f)$ must normalise its characteristic subgroup $Z(f)$, so $N_G(Q(f))\le N_G(Z(f))=QA_0(f)$. Now $A_0(f)\le F(f)\le N_G(Q(f))$, and \S3.4.1(e) shows that $N_Q(Q(f))=\bigcup_{\alpha}(\alpha,0)Z$ with the union over all $\alpha\in{\mathbb F}(f)$, so $N_G(Q(f))$ has order $|Z|.2^f.|A_0(f)|=2^{e+f}(2^f-1)$.

\smallskip

\noindent(2) Clearly $N_G(F(f))\le N_G(Q(f))\le QA_0(f)$, and $A_0(f)\le N_G(F(f))$. Since $Z$ acts semi-regularly on $\Omega\setminus\{\infty\}$, it acts semi-regularly by conjugation on the subgroups of $F$ in the conjugacy class $\mathcal{A}_0$, so
$N_G(F(f))\cap Z=Z(f)$. Hence, using the proof of part~(3), we see that $N_G(F(f))\le F\cap G(f)=F(f)$. Thus $F(f)$ is self-normalising.

\medskip

\noindent(5, 6) See \S 2.3.2.

\smallskip

\noindent(7, 8) For $i=0$ and $2$ the subgroups in ${\mathcal B}_i(1)$ are cyclic groups of orders $2$ and $4$ respectively, so they are contained in Sylow $2$-subgroups of $G$. There are $q^2+1$ Sylow $2$-subgroups, each conjugate to $Q$ and containing $q-1$ subgroups of order $2$, and containing $(q^2-q)/2$ of order $4$. Since distinct Sylow $2$-subgroups have trivial intersection, there are $(q^2+1)(q-1)$ and $(q^2+1)(q^2-q)/2$ such subgroups in $G$. In each case such subgroups are all conjugate, so their normalisers have order $q^2$ and $2q$. \hfill$\square$

\section{Calculating values of $\mu_G$}

We can now complete the proof of Theorem~\ref{mainthm} by calculating $\mu_G(H)$ for each subgroup $H\in{\mathcal S}$. In order to use equation~(\ref{recursive}) for this (see \S 1), we first need to know, for each pair of subgroups $H, K\in{\mathcal S}$, the number $N(H;K)$ of conjugates in $G$ of $K$ containing $H$. If $M(H;K)$ denotes the number of conjugates in $G$ of $H$ contained in $K$, then a simple double counting argument gives
\begin{equation}\label{MN}
M(H;K)M(K;G) = M(H;G)N(H;K)
\end{equation}
for all $H, K\in{\mathcal S}$. This allows $N(H;K)$ to be determined from the values of the function $M$. Now $M(H;G)=|{\mathcal H}|$ and $M(K;G)=|{\mathcal K}|$, where $\mathcal H$ and $\mathcal K$ are the conjugacy classes of subgroups of $G$ containing $H$ and $K$, so these values are given by Theorem~\ref{conjsize}. The values of $M(H;K)$ for $K\ne G$ can be found by using arguments similar to those used in proving Theorem~\ref{conjsize}, so details are omitted.

\begin{table}[ht]\label{Ntable}

\centering

\begin{tabular}{| p{3cm} | p{1cm} | p{1cm} | p{1cm} | p{1cm} | p{1cm} | p{1cm} |}
\hline
$~$ & $G(h)$ & $F(h)$ & $Q(h)$ & $Z(h)$ & $B_0(h)$ & $A_0(h)$ \\
\hline\hline
$G(f)$ & $1$ & & & & & \\
\hline
$F(f)$ & $1$ & $1$ & & & & \\
\hline
$Q(f)$ & $2^{e-h}$ & $2^{e-h}$ & $1$ & & & \\
\hline
$Z(f)$ & $2^{2(e-h)}$ & $2^{2(e-h)}$ & $2^{e-h}$ & $1$ & & \\
\hline
$B_0(f)$ & $1$ & & & & $1$ & \\ 
\hline
$A_0(f)$ & $\frac{(2^e-1)}{(2^h-1)}$ & $\frac{2(2^e-1)}{(2^h-1)}$ & & & $\frac{(2^e-1)}{(2^h-1)}$ & $1$ \\ [0.8ex]
\hline
$B_1(f), f\geqslant 1$ & $1$ & & & & & \\
\hline
$B_2(f)$ & $1$ & & & & & \\
\hline
$A_1(f), f\geqslant 1$ & $\frac{a_1(e)}{a_1(h)}$ & & & & & \\ [0.8ex]
\hline
$A_2(f)$ & $\frac{a_2(e)}{a_2(h)}$ & $~$ & $~$ & $~$ & $~$ & $~$ \\ [0.8ex]
\hline
$B_2(1)\cong C_4$ & $2^{e-h}$ & $2^{e-h}$ & $1$ & & & \\
\hline
$B_0(1)\cong C_2$ & $2^{2(e-h)}$ & $2^{2(e-h)}$ & $2^{e-h}$ & $1$ & $2^{2e-1}$ & \\
\hline

\end{tabular}
\label{t1}
\end{table}

\begin{table}[ht]
\centering
\begin{tabular}{| p{3cm} | p{1cm} | p{1cm} | p{1cm} | p{1cm} | p{1cm} | p{1cm} |}
\hline
$~$ & $B_1(h)$ & $B_2(h)$ & $A_1(h)$ & $A_2(h)$ & $B_2(1)$ & $B_0(1)$ \\
\hline\hline
$B_1(f), f\geqslant 1$ & $1$ & & & & & \\
\hline
$B_2(f)$ & & $1$ & & & & \\
\hline
$A_1(f), f\geqslant 1$ & $\frac{a_1(e)}{a_1(h)}$ & $~$ & $1$ & $~$ & $~$ & $~$\\[0.8ex]
\hline
$A_2(f)$ & & $\frac{a_2(e)}{a_2(h)}$ & & $1$ & & \\[0.8ex]
\hline
$B_2(1)\cong C_4$ & $2^{e-1}$ & $2^{e-1}$ & $~$ & $~$ & $1$ & $~$ \\
\hline
$B_0(1)\cong C_2$ & $2^{2(e-1)}$ & $2^{2(e-1)}$ & $~$ & $~$ & $2^{e-1}$ & $1$ \\
\hline

\end{tabular}

\caption{Values of $N(H;K)$ where $H, K\in{\mathcal S}$.}

\end{table}

The non-zero values of $N(H;K)$ resulting from~(\ref{MN}) are given in Table~2, where the rows and columns are indexed by the subgroups $H$ and $K$ respectively; the row corresponding to the identity subgroup $H=A_0(1)=A_2(1)$ is omitted since in this case $N(H;K)=|{\mathcal K}|$, given by Theorem~\ref{conjsize} for all $K\in{\mathcal S}$. The table is split into two parts, the second part giving further entries for the last six rows of the first part. We assume that $f$ divides $h$ and that $f>1$ unless otherwise stated. Thus $G(1)$ is represented by its conjugate $B_1(1)$, while $F(1)$ and $Q(1)$ are represented by $B_2(1)$, and $Z(1)$ by $B_0(1)$ (see \S2.3.3 and the comments in \S2.4}).

Given Table~2, one can systematically use equation~(\ref{recursive}) to calculate $\mu_G(H)$ for each $H\in{\mathcal S}$, starting with $H=G(f)$ in the first row, and working downwards through the table. For instance, if $H=G(f)$ then the subgroups $K\in{\mathcal S}$ with $N(H;K)\ne 0$ are those of the form $K=G(h)$ where $f\mid h\mid e$; under inclusion, these form a lattice isomorphic to the lattice $\Lambda(e/f)$ of all divisors $h/f$ of $e/f$, with $\mu_G(K)=1$ when $h=e$, so we find that $\mu_G(H)=\mu(e/f)$, as in Table~1. Next, if $H=F(f)$ we consider the subgroups $K=G(h)$ and $F(h)$ where $f\mid h\mid e$; these form a lattice isomorphic to $\Lambda(2e/f)$ since $e$ is odd, giving $\mu_G(H)=\mu(2e/f)=-\mu(e/f)$. Similar arguments show that if $f>1$ then $\mu_G(B_i(f))=-\mu(e/f)$ for $i=0, 1, 2$, and $\mu_G(Q(f))=\mu_G(Z(f))=\mu_G(A_i(f))=0$ for $i=1, 2$. This process continues until $\mu_G(H)$ is evaluated for all $H\in{\mathcal S}$. The method is essentially the same as that described fully in~\cite[\S4]{DowJLMS} for the groups $G=L_2(2^e)$, so the remaining details are omitted.

Since $\mu_G(H)=0$ whenever $H=Q(f)$, $Z(f)$, $A_1(f)$ or $A_2(f)$ for any $f>1$, we may disregard these subgroups, and let $\mathcal T$ denote the remaining set of subgroups $H\in{\mathcal S}$, namely those of the form
\begin{equation}\label{setT}
G(f),\; F(f),\; B_i(f)\; (i=0, 1, 2),\; A_0(f),\; B_2(1),\; B_0(1),\; A_0(1),
\end{equation}
where $1<f\mid e$. This is a set of representatives for the conjugacy classes in Table~1; by the construction of $\mathcal T$, every subgroup $H$ of $G$ with $\mu_G(H)\ne 0$ must belong to one of these classes. This fact, together with the values of $|H|$, $|N_G(H)|$ and $\mu_G(H)$ determined earlier, justifies the entries in Table~1 and in particular proves Theorem~\ref{mainthm}. \hfill$\square$

\medskip

Each conjugacy class in Table~1 contains $|G|/|N_G(H)|$ subgroups, and $|{\rm Aut}\,G|=e|G|$ by \S\ref{basic}(2), so equation~(\ref{NGammaG}) takes the form
\[
n_{\Gamma}(G)
=\frac{1}{e}\sum_{H\in{\mathcal T}}
\frac{\mu_G(H)|{\rm Hom}(\Gamma,H)|}{|N_G(H)|},
\]
as in equation~(\ref{NewNGammaG}). Table~1 gives the values of  $\mu_G(H)$ and $|N_G(H)|$, so in order to apply this equation to a particular group $\Gamma$ one needs only to count the homomorphisms $\Gamma\to H$ for each $H\in{\mathcal T}$.
We will illustrate this in the case of certain Hecke groups $\Gamma$ in Section~\ref{Hecke}.

\section{Application to Hecke groups}\label{Hecke}

In~\cite{DowJLMS}, the first author used M\"obius inversion to enumerate the normal subgroups of the modular group $\Gamma=PSL_2({\mathbb Z})\cong C_2*C_3$ with quotient group isomorphic to $L_2(q)$ for any given prime power $q$. There are no normal subgroups of $\Gamma$ with quotient group $G=Sz(q)$, since Suzuki groups have no elements of order $3$, but instead one can apply the same method to Hecke groups other than $\Gamma$.

For each integer $k\ge 3$ the Hecke group $H_k$, introduced by Hecke in connection with Dirichlet series~\cite{Hec}, is the subgroup of $PSL_2({\mathbb R})$ generated by the M\"obius transformations
\[X: z\mapsto z+\lambda_k  \quad {\rm and} \quad Y: z\mapsto \frac{-1}{z} \]
of the upper half plane, where $\lambda_k=2\cos(\pi/k)$. This group is a free product of cyclic groups of orders $2$ and $k$ generated by $Y$ and $XY$. In particular, $H_3$ is the modular group.

The torsion-free normal subgroups of $H_k$ with a given finite quotient group $G$ correspond bijectively to the orbits of ${\rm Aut}\,G$ on generating pairs for $G$ of orders $2$ and $k$. By equation~(\ref{NewNGammaG}) the number of these orbits is
\begin{equation}\label{Hk}
\frac{1}{|{\rm Aut}\,G|}\sum_{H\le G}\mu_G(H)|H|_2|H|_k\;
=\frac{1}{e}\sum_{H\in{\mathcal T}}
\frac{\mu_G(H)|H|_2|H|_k}{|N_G(H)|},
\end{equation}
where $|H|_n$ is the number of elements of order $n$ in each subgroup $H$ of $G$. Table~3 gives $|H|_n$ for subgroups $H\in{\mathcal T}$ in the cases $n=2, 4$ and $5$; here $q=2^e$ and $1<f\mid e$. 

\begin{table}[ht]
\centering
\begin{tabular}{| p{2.1cm} | p{2.7cm} | p{3.5cm} | p{3cm} |}
\hline
$H$  & $|H|_2$ & $|H|_4$ & $|H|_5$ \\
\hline\hline
$G(f)$ &   $(2^f-1)(2^{2f}+1)$ & $2^f(2^{2f}+1)(2^f-1)$ & $2^{2f}(2^f-1)a_2(f)$\\
\hline
$F(f)$  & $2^f-1$ & $2^f(2^f-1)$ & $0$ \\
\hline
$B_0(f)$  & $2^f-1$ & $0$ & $0$ \\
\hline
$A_0(f)$  & $0$ & $0$ & $0$ \\ 
\hline
$B_1(f)$  & $a_1(f)$ & $2a_1(f)$ & $4$ \\
\hline
$B_2(f)$  & $a_2(f)$ & $2a_2(f)$ & $0$ \\
\hline
$B_2(1)$ & $1$ & $2$ & $0$ \\
\hline
$B_0(1)$ & $1$ & $0$ & $0$ \\
\hline
$A_0(1)$ & $0$ & $0$ & $0$ \\
\hline

\end{tabular}
\caption{Values of $|H|_n$ for $n=2, 4$ and $5$.}
\label{table:|H|k}
\end{table}

With this information, together with Table~1, one finds from equation~(\ref{Hk}) that the numbers of (necessarily torsion-free) normal subgroups of $H_4$ and of $H_5$ with quotient group $Sz(q)$ are respectively
\begin{equation}\label{H4}
\frac{1}{e}\sum_{f|e}\mu\left(\frac{e}{f}\right)2^f(2^f-2)
\end{equation}
and
\begin{equation}\label{H5}
\frac{1}{e}\sum_{f|e}\mu\left(\frac{e}{f}\right)(2^f-1)a_2(f),
\end{equation}
where $a_2(f)=|A_2(f)|=2^f-\chi(f)2^{(f+1)/2}+1$ (see \S\ref{defn}).

One can apply similar arguments for odd $k>5$. For instance, $G$ contains elements of order $7$ if and only if $e$ is divisible by $3$, in which case they form three conjugacy classes, represented by elements of $A_0$. It follows that the number of normal subgroups of $H_7$ with quotient group $Sz(q)$ is
\begin{equation}\label{H7}
\frac{3}{e}\sum_{3|f|e}\mu\left(\frac{e}{f}\right)(2^{2f}-2).
\end{equation}

These enumerations have applications in other areas, such as topological graph theory.  For instance, formulae~(\ref{H4}), (\ref{H5}) and~(\ref{H7}) give the numbers of orientably regular $k$-valent maps with automorphism group $Sz(q)$ for $k=4, 5$ and $7$  (see~\cite{DJ13} for these and other examples).

\end{document}